\documentclass{article}

\usepackage{xr,refcount}
\usepackage{amsmath}
\usepackage{amssymb}
\usepackage{amsthm}
\usepackage{cite}
\usepackage{hyperref}
\usepackage{enumerate}
\usepackage{amsrefs}

\newcommand{\ZZ}[0]{\mathbb{Z}}

\newcommand{\EE}[0]{\mathbb{E}}

\newcommand{\skp}[0]{\mathsf s_{kp}}

\newtheorem{theorem}{Theorem}

\newtheorem{conj}{Conjecture}
\newtheorem{definition}{Definition}

\theoremstyle{remark}

\usepackage{listings}
\usepackage{color}

\definecolor{dkgreen}{rgb}{0,0.6,0}
\definecolor{gray}{rgb}{0.5,0.5,0.5}
\definecolor{mauve}{rgb}{0.58,0,0.82}

\title{Exponential Lower Bounds on the Generalized Erd{\H o}s-Ginzburg-Ziv Constant}
\date{September 1, 2017}
\author{Jared Bitz, Sarah Griffith, and Xiaoyu He}

\begin{document}
\maketitle

\begin{abstract}
For a finite abelian group $G$, the generalized Erd\H{o}s--Ginzburg--Ziv constant
$\mathsf s_{k}(G)$ is the smallest $m$ such that a sequence of $m$ elements in $G$
always contains a $k$-element subsequence which sums to zero. If $n = \exp(G)$ is the exponent of $G$, the previously best known bounds for $\mathsf s_{kn}(C_n^r)$ were linear in $n$ and $r$ when $k\ge 2$.
Via a probabilistic argument, we produce the exponential lower bound 
\[
\mathsf s_{2n}(C_n^r) > \frac{n}{2}[1.25 + o(1)]^r 
\]
for $n > 0$.
For the general case, we show 
\[
\mathsf s_{kn}(C_n^r) > \frac{kn}{4}\Big(1+\frac{1}{ek+1} + o(1)\Big)^r.
\]
\end{abstract}

\section{Introduction}
\label{sec:intro}
In 1961, Erd\H{os}, Ginzburg, and Ziv \cite{Erdos1961} proved that among any $2n - 1$ integers, some $n$ of them sum to a multiple of $n$.
Equivalently, among every sequence (with repetition) of $2n - 1$ elements in $C_n := \ZZ/n\ZZ$, some $n$ sum to zero (for brevity, given a sequence, we call a subsequence of length $n$ an \textit{$n$-subsequence}).

This result led naturally to the study of restricted-length zero-sum subsequences in finite abelian groups.
Recall that the exponent $\exp(G)$ of a finite group is the largest order among its elements.

\begin{definition}
If $G$ is a finite abelian group and $n=\exp(G)$, the $k$-th generalized EGZ constant of $G$, denoted $\mathsf{s}_{kn}(G)$, is the smallest $m$ for which any sequence of $m$ elements of $G$ contains a zero-sum $kn$-subsequence.
\end{definition} 


The EGZ problem, especially in the ``smallest'' case of $k = 1$, has proven to be surprisingly difficult.
The only $r$ for which $s_n(C_n ^r)$ has been exactly determined are $r=1$ and $r=2$.
The fact that $\mathsf{s}_{n}(C_n) = 2n - 1$ is the original Erd\H os-Ginzburg-Ziv Theorem \cite{Erdos1961}, while $\mathsf{s}_{n}(C_n^2) = 4n - 3$ is the famous Kemnitz conjecture and was only recently settled by Reiher \cite{Reiher2007}. Reiher's theorem builds on the polynomial method of R\'onyai \cite{Ronyai2000}, which has proved fruitful for providing upper bounds to $\mathsf{s}_{kn}(C_n^r)$ in higher ranks $r\ge 3$, especially when $k$ is somewhat large compared to $r$.

Major progress on the hardest case $k=1$ was recently made by Naslund  \cite{Naslund2017}, obtaining an exponential improvement on the upper bounds on $\mathsf{s}_p(C_p^r)$ for general $r$ by using the
``multi-slice-rank'' method.
This technique generalizes that of Ellenberg and Gijswijt \cite{Ellenberg2016} and Croot, Lev, and Pach \cite{Croot2017} on the closely related cap-set problem.
Specifically, Naslund shows that
\begin{equation*}
    \mathsf{s}_p(C_p^r) < 3p!(2p-1)(J(p)p)^r,
\end{equation*}
where $J(p)$ is a constant, depending on $p$, which lies between .841 and .918 and decreases as $p$ grows.

We turn our attention to $\mathsf{s}_{kn}(C_n^r)$ for higher values of $k$. Gao and Thangadurai \cite{Gao2006} proved that $5p + \frac{p-1}{2} -3 \leq 
\mathsf{s}_{2p}(C_p^3) \leq 6p - 3$, and that $\skp(C_p^3) = (k+3)p - 3$ for $k \geq 4$.
Kubertin \cite{kubertin2005} extended the latter result to the $k = 3$ case, and also showed that
$s_{4p}(C_p^4) = 8p - 4$.
These results contribute to the following conjecture, made in a more general form by Gao and Thangadurai \cite{Gao2006}.

\begin{conj}
For all $k \geq r$, and sufficiently large $n$, $\mathsf{s}_{kn}(C_n^r) = (k + r)n - r$.
\end{conj}

Extending the work of R\'onyai\cite{Ronyai2000} and Kubertin\cite{kubertin2005}, it was shown by He \cite{He2015} that the conjecture holds for prime $n$ when $k \geq p + r$ and $2p \geq 7r - 3$. For a survey of related zero-sum problems and their generalizations, see \cite{GaoSurvey}.

When $n$ is an odd prime and $k=1$, Harborth \cite{Harborth} gives the elementary lower bound 
$\mathsf{s}_p(C_p^r) > 2^r(p-1)$, but no similar bound exponential in $r$ was known for $k \ge 2$. In this paper, we provide lower bound constructions for $\mathsf s_{kn}(C_n^r)$ when $k\ge 2$ which are effective when $k$ is much smaller than $r$.

\begin{theorem}
    \label{thm:theorem1}
    The generalized EGZ constant $\mathsf{s}_{2n}(C_n^r)$ satisfies the bound
    \begin{equation*}
        \mathsf{s}_{2n}(C_n^r) > \frac{n}{2}\left[\frac{5}{4} + o(1)\right]^r.
    \end{equation*}
\end{theorem}
That is, there is a sequence of this length in $C_n^r$ that
contains no zero-sum $2n$-subsequence.
This result is the first exponential lower bound on $\mathsf{s}_{2n}(C_n^r)$ for an arbitrary $r$.

Generalizing this to arbitrary values of $k$, we also prove
\begin{theorem}
	\label{thm:theorem2}
    Let $k > 2$. 
    The generalized EGZ constant $\mathsf{s}_{kn}(C_n^r)$ satisfies the bound
    \begin{equation*}
    	\mathsf{s}_{kn}(C_n^r) > \frac{kn}{4}\left[1 + \frac{1}{ek+1} + 
        o(1)\right]^r.
    \end{equation*}
\end{theorem}
In both results, the $o(1)$ term goes to zero as $n$ grows.

\section{Lower Bounds on $\mathsf{s}_{kn}(C_n^r)$}
\label{sec:lowerbounds}
We first give the proof of \autoref{thm:theorem1}.
\begin{proof}[Proof of \autoref{thm:theorem1}]
    Let
    \begin{equation*}
        N = \frac{n}{2}A^r
    \end{equation*}
    for a value of $A$ which will be specified later.
    Choose a sequence $X$ of $N$ random vectors in $\{0, 1\}^r$ as follows.
    For each $v = (v_1, \dots, v_r)$ that is a term of $X$, 
    let each $v_i = 1$ with probability $q$ and $v_i = 0$ with
    probability $1-q$, with each $v_i$ chosen independently.
    Let $Z$ be the number of zero-sum length $2n$-subsequences in $X$.
    We will produce an $A$ such that $\EE[Z] < 1$, so that there must be some possible $X$ with
    no such subsequences.

    Consider some arbitrary $2n$-subsequence $Y$ of $X$.
    For $Y$ to be zero-sum, each of the $r$ coordinates must sum to 0 mod $n$,
    and so contain exactly 0, $n$, or $2n$ ones.
    For any coordinate $i \leq r$, let $P_0$ be the probability that coordinate
    $i$ contains 0 ones, $P_n$ the probability that it contains $n$ ones, and
    $P_{2n}$ the probability that it contains $2n$ ones (clearly, this is not dependent on $i$).
    We have
    \begin{equation*}
        \begin{split}
            P_0 &= (1-q)^{2n} \\
            P_n &= {2n \choose n}q^n(1-q)^n\\
            P_{2n} &= q^{2n},
        \end{split}
    \end{equation*}
    the values from a standard binomial distribution with probability of 
    success $q$.

    Now, if $Q$ is the probability that coordinate $i$ sums to zero, we have
    \begin{equation*}
        Q = P_0 + P_n + P_{2n}.
    \end{equation*}
    We wish to minimize $Q$. We proceed by choosing $q$ so that $P_0$ grows slowly while still dominating. Calculations indicate that we ought to allow $q$ to equal $4/5$.
    Using the bound ${2n \choose n} \leq 4^n/\sqrt{3n + 1}$, we then have
    \begin{equation*}
        Q \leq \left(\frac{4}{5}\right)^{2n} + \frac{1}{\sqrt{3n+1}}\left(\frac{4}{5}\right)^{2n}
        + \left(\frac{1}{5}\right)^{2n} < 
        \left(1 + \frac{1}{\sqrt{n}}\right)\left(\frac{4}{5}\right)^{2n}.
    \end{equation*}
    Since $Q$ is the probability that any one coordinate in $Y$ sums to zero,
    the probability that $Y$ as a whole is zero-sum is $Q^r$.
    Since each $2n$-subsequence of $X$ is zero-sum with equal probability, we
    have
    \begin{equation*}
        \begin{split}
            \EE[Z] &= {N \choose 2n}Q^r \\
                   &< \left(\frac{2N}{n}\right)^{2n}
                   		\left(1 + \frac{1}{\sqrt{n}}\right)^r\left(\frac{4}{5}\right)^{2nr} \\
                   &< A^{2nr}\left(1 + \frac{1}{\sqrt{n}}\right)^r\left(\frac{4}{5}\right)^{2nr}. \\
        \end{split}
    \end{equation*}
    To obtain $\EE[Z] < 1$, it is sufficient to have
    \begin{equation*}
    \begin{split}
    	A^{2nr}\left(1 + \frac{1}{\sqrt{n}}\right)^r\left(\frac{4}{5}\right)^{2nr} &< 1 \\
        A < \frac{5}{4}\left(1 + \frac{1}{\sqrt{n}}\right)^{-\frac{1}{2n}},
    \end{split}
    \end{equation*}
    So, we can take $A = \frac{5}{4} + o(1)$.
\end{proof}

It is possible to improve \autoref{thm:theorem1} by a subexponential factor using the Lovasz Local Lemma, but the exponential constant $5/4$ appears to be the natural limit of the method.

The proof of \autoref{thm:theorem2} is similar.

\begin{proof}[Proof of \autoref{thm:theorem2}]
If we repeat the above construction in the general case, 
we first must calculate the probability
that a given coordinate $i \leq r$ sums to zero.
We now have $k+1$ component probabilities, $P_0$ through $P_{kn}$.
By the same logic as above, we see that
\begin{equation*}
	P_{in} = {kn \choose in} q^{in}(1-q)^{(k-i)n}.
\end{equation*}
We again wish to have $P_0$ dominate while growing slowly. This time calculations suggest allowing $q$ to equal $1/(ek + 1)$. Note that when $P_0$ dominates we have the bound $Q < (k + 1)(1-q)^{kn}$.
Then, if $N = \frac{kn}{4}A^r$, we have
\begin{equation*}
    \begin{split}
        \EE[Z] = {N \choose k}Q^r &< \left(\frac{4N}{kn}\right)^{kn}
        (k+1)^r(1-q)^{knr}.
    \end{split}
\end{equation*}
If we want $\EE[Z] < 1$, we must have
\begin{equation*}
        \left(\frac{4N}{kn}\right)^{kn}(k+1)^r(1-q)^{knr} < 1
\end{equation*}
and thus
\begin{equation*}
    \begin{split}
        A &< \frac{1}{(k+1)^{1/kn}(1-q)} \\
          &< \frac{1}{(k+1)^{1/kn}\left(1 - \frac{1}{ek + 1}\right)}.
    \end{split}
\end{equation*}
So, we have $A = 1 + \frac{1}{ek + 1} + o(1)$, as desired.
\end{proof}
When $q$ is determined as above, the permissible values of $A$
approach one as $k$ and $n$ grow.
Therefore, the result gives us little if we seek a bound on $\mathsf s_{kn}(C_n^r)$
for completely general $k$ and $n$.
However, if we fix a specific value of $k$ (for example, if we are interested in $\mathsf{s}_{3n}$), we can
still attempt to optimize $q$ in order to achieve an exponential lower bound - 
it just will not be as effective for large $k$.

\section{Acknowledgements}
The authors would like to thank George Schaeffer and the Stanford Undergraduate
Research in Mathematics program for support during the development of these results.
We indebted to Jacob Fox for valuable guidance on the probabilistic method,
and to Christian Reiher for general insights regarding the Erd{\H o}s-Ginzburg-Ziv problem. We would also like to thank Jesse Geneson for the helpful comments he provided us with.

\bibliographystyle{plain}
\bibliography{bibliography}

\begin{bibdiv}
\begin{biblist}

\bib{Croot2017}{article}{
      author={Croot, E.},
      author={Lev, V.},
      author={Pach, P.},
       title={{Progression-free sets in $\mathbb{Z}_4^n$ are exponentially
  small}},
        date={2017},
        ISSN={0003-486X},
     journal={Annals of Mathematics},
      volume={185},
       pages={331\ndash 337},
  url={http://arxiv.org/abs/1605.01506{\%}5Cnhttp://annals.math.princeton.edu/2017/185-1/p07},
}

\bib{Ellenberg2016}{article}{
      author={Ellenberg, J.},
      author={Gijswijt, D.},
       title={{On large subsets of $\mathbb{F}_q^n$ with no three-term
  arithmetic progressions}},
        date={2016},
        ISSN={0003486X},
     journal={Annals of Mathematics},
      number={185},
       pages={1\ndash 4},
}

\bib{Erdos1961}{article}{
      author={Erd\H{o}s, P.},
      author={Ginzburg, A.},
      author={Ziv, A.},
       title={{Theorem in the additive number theory}},
        date={1961},
     journal={Bull. Research Council Israel},
      volume={10},
       pages={41\ndash 43},
}

\bib{GaoSurvey}{article}{
      author={Gao, W.},
      author={Geroldinger, A.},
       title={{Zero-sum problems in finite abelian groups: A survey}},
        date={2006},
        ISSN={07230869},
     journal={Expositiones Mathematicae},
      volume={24},
       pages={337\ndash 369},
}

\bib{Gao2006}{article}{
      author={Gao, W.},
      author={Thangadurai, R.},
       title={{On zero-sum sequences of prescribed length}},
        date={2006},
     journal={Aequationes Mathematicae},
      volume={72},
       pages={201\ndash 212},
}

\bib{Harborth}{article}{
      author={Harborth, H.},
       title={{Ein Extremalproblem fur Gitterpunkte}},
        date={1973},
     journal={J. Reine Angew. Math.},
      number={262},
       pages={356\ndash 360},
}

\bib{He2015}{article}{
      author={He, X.},
       title={Zero-sum subsequences of length kq over finite abelian p-groups},
        date={2015},
     journal={Discrete Mathematics},
      number={339},
       pages={399\ndash 407},
}

\bib{kubertin2005}{article}{
      author={Kubertin, S.},
       title={Zero-sums of length $kq$ in $\mathbb{Z}_q^d$},
        date={2005},
     journal={Acta Arithmetica},
      volume={116},
       pages={145\ndash 152},
}

\bib{Naslund2017}{article}{
      author={Naslund, E.},
       title={{Exponential bounds for the Erd\H{o}s-Ginzburg-Ziv constant}},
        date={2017},
     journal={preprint},
      eprint={arXiv:1701.04942v1},
}

\bib{Reiher2007}{article}{
      author={Reiher, C.},
       title={{On Kemnitz' conjecture concerning lattice-points in the plane}},
        date={2007},
     journal={Ramanujan Journal},
      volume={13},
       pages={333\ndash 337},
}

\bib{Ronyai2000}{article}{
      author={R\'onyai, L.},
       title={{On a Conjecture of Kemnitz}},
        date={2000},
        ISSN={0209-9683},
     journal={Combinatorica},
      volume={20},
       pages={569\ndash 573},
         url={http://link.springer.com/10.1007/s004930070008},
}

\end{biblist}
\end{bibdiv}

\end{document}